\documentclass{amsart}

\usepackage{amsmath}
\usepackage{amsfonts}
\usepackage{amscd}
\usepackage{amssymb}
\input xypic

\newtheorem{defn}{Definition}[section]
\newtheorem{thm}[defn]{Theorem}
\newtheorem{prop}[defn]{Proposition}
\newtheorem{lemma}[defn]{Lemma}

\newtheorem{cor}[defn]{Corollary}

\newtheorem{schol}[defn]{Scholium}

\newcommand{\lm}{\ensuremath{\longrightarrow}}

\newcommand{\eps}{\varepsilon}

\DeclareMathOperator{\Hom}{\mbox{Hom}}
\DeclareMathOperator{\shom}{\ensuremath{\mathcal{H}\mathit{om}}}
\DeclareMathOperator{\Rshom}{\ensuremath{R\mathcal{H}\mathit{om}}}

\DeclareMathOperator{\send}{\ensuremath{\mathcal{E}\!\mathit{nd}}}
\DeclareMathOperator{\sext}{\ensuremath{\mathcal{E}\!\mathit{xt}}}
\DeclareMathOperator{\End}{\mbox{End}}
\DeclareMathOperator{\Aut}{\mbox{Aut}}

\DeclareMathOperator{\Ext}{\mbox{Ext}}

\DeclareMathOperator{\id}{\mbox{id}\,}

\DeclareMathOperator{\im}{\mbox{im}\,}

\DeclareMathOperator{\Pic}{\mbox{Pic}\,}

\DeclareMathOperator{\Div}{\mbox{Div}\,}
\DeclareMathOperator{\rk}{\mbox{rank}\,}

\DeclareMathOperator{\s}{\sigma}
\DeclareMathOperator{\z}{\zeta}
\DeclareMathOperator{\w}{\omega}

\DeclareMathOperator{\PP}{\mathbb{P}}
\DeclareMathOperator{\Z}{\mathbb{Z}}

\DeclareMathOperator{\calo}{\mathcal{O}}

\DeclareMathOperator{\oy}{\mathcal{O}_{Y}}

\title{Moduli of Bundles on Exotic del Pezzo Orders}

\author{Daniel Chan} 
\author{Rajesh S. Kulkarni}
\address{School of Mathematics, University of New South Wales, Sydney, NSW 2052, Australia}
\email{danielc@unsw.edu.au} 
\address{Department of Mathematics, Wells Hall, Michigan State University, East Lansing, MI 48824, USA}
\email{kulkarni@math.msu.edu}

\begin{document}
\begin{abstract}
We study bundles on rank 4 maximal orders on $\PP^2$ ramified on a smooth plane quartic. We compute the possible Chern classes for line bundles. Our main result is that the moduli space of line bundles with minimal second Chern class is either a point or a smooth genus two curve.
\end{abstract}
\maketitle
Throughout, all objects and maps are assumed to be defined over some algebraically closed base field $k$ of characteristic 0. Modules will by default be left modules unless otherwise noted.

\section{Introduction}  \label{sintro}

The Picard variety and more generally, moduli spaces of vector bundles form important invariants for projective varieties. Simpson and others (see \cite{Simp},\cite{HoSt},\cite{Lieb}) have noted that much of the general theory of moduli spaces of sheaves extends to finite sheaves of non-commutative algebras over a projective variety. However, little work has been done in studying specific moduli spaces, mainly because interesting examples of such sheaves of algebras are hard to describe explicitly.

In \cite{C}, examples of orders on projective varieties were constructed using a non-commutative analogue of the cyclic covering trick, leading to what is called a non-commutative cyclic cover (the definition of which will be given in the section~\ref{snccovers}). One particular example of interest is a construction of a degree two maximal order $A$ on $\PP^2$ ramified on a smooth quartic $D$. It is one of the exotic del Pezzo orders. The main objective of this paper is to study ``line bundles'' on this order. It is important to note that these line bundles are one-sided $A$-modules so cannot be tensored together to form a group structure. In particular, the moduli space of these is not naturally a group scheme.

To explain the results, recall that $A$ can be written as a non-commutative cyclic cover of the del Pezzo surface $Y$ where $Y$ is the double cover of $\PP^2$ ramified on $D$ (so $Y$ and $A$ have the same ramification over $\PP^2$). In particular, $\oy$ is a maximal commutative subalgebra of $A$, though $A$ is not a sheaf of algebras on $Y$ in the usual sense. Hence $A$-modules are also sheaves on $Y$, so we may consider their Chern classes as sheaves on $Y$. An $A$-line bundle is an $A$-module of rank 2 on $Y$.  We classify the possible Chern classes of $A$-line bundles. We show for fixed $c_1$, there is a lower bound on $c_2$, but the second Chern class is otherwise unrestricted. We were rather surprised to find out that the coarse moduli scheme of $A$-line bundles is computable in the case when $c_2$ is minimal. There are two possibilities, the moduli scheme is either a point, or is a double cover of $\PP^1$ ramified at 6 points.

We hope this paper is of interest to, and readable by, both algebraic geometers with minimal knowledge of non-commutative algebra, and ring theorists familiar with some algebraic geometry. To this end, we review non-commutative cyclic covers in section~\ref{snccovers}. We hope this section will in particular, contain most of the non-commutative algebra an algebraic geometer needs to read the paper. Our approach to studying $A$-line bundles is by examining rank 2 bundles on $Y$ with extra structure. This is outlined in section~\ref{slineb}. The construction of the exotic del Pezzo orders on $\PP^2$ ramified on a smooth quartic is recalled in section~\ref{sdPcase} where we also verify smoothness of the moduli scheme of line bundles. The last two sections compute the Chern classes and moduli schemes as alluded to above.

\section{Review of Non-commutative Cyclic Covers} \label{snccovers}
In this section, we review the definition of non-commutative cyclic covers and give some of the basic properties that we will need. For more details see \cite{C} and \cite{LVV}.

Let $Y$ be a quasi-projective variety. The most elegant definition of the non-commutative cyclic cover involves van den Bergh's notion of an $\oy$-bimodule (see \cite[Section~2]{VdB}). We will give this version of the definition first. For the sake of the reader familiar with only commutative algebraic geometry, we will also give another description of non-commutative cyclic covers which bypasses this theory.

Recall that Van den Bergh constructs in \cite[Section~2]{VdB} a monoidal category of quasi-coherent $\oy$-bimodules. The invertible objects have the form $L_{\s}$ where $L \in \Pic Y$ and $\s \in \Aut Y$. One may think of this intuitively, as the $\oy$-module $L$ where the right module structure is skewed through by $\s$ so that $\ _{\oy} L \simeq L, L_{\oy} \simeq \s_* L$.

Suppose now that $\s$ generates a finite cyclic group $G$ of order dividing some integer $e$, and that there is a non-zero map of bimodules $\phi: L_{\s}^{\otimes e} \lm \oy$ which we view as a relation. As in the commutative cyclic covering trick, we wish to put an algebra structure on
$$ A(Y;L_{\s},\phi) := \oy \oplus L_{\s} \oplus L_{\s}^{\otimes 2} \oplus \ldots \oplus L_{\s}^{\otimes(e-1)} $$
using $\phi$ to account for the non-obvious multiplication maps. It turns out that this is possible if $\phi$ satisfies the {\em overlap condition}, namely, if the two maps $1\otimes \phi, \phi \otimes 1: L_{\s}^{\otimes e+1} \lm L_{\s}$ are equal. In this case, we refer to $A(Y;L_{\s},\phi)$ as a {\em cyclic algebra} which gives an $e$-fold {\em non-commutative cyclic cover} of $Y$. The multiplication is given by the natural isomorphism $L_{\s}^{\otimes i} \otimes L_{\s}^{\otimes j} \lm L_{\s}^{\otimes i+j}$ when $i+j < e$ and 
$$L_{\s}^{\otimes i} \otimes L_{\s}^{\otimes j} \lm L_{\s}^{\otimes i+j} \xrightarrow{1 \otimes \phi\otimes 1} L_{\s}^{\otimes i+j-e}$$ 
when $i+j \geq e$. The overlap condition implies that the definition is independent of which $e$ consecutive tensor factors of $L_{\s}$ you apply $\phi$ to in the above expression. 

We can re-interpret this algebra without using bimodules as follows. We need to know that for invertible bimodules $L_{\s},M_{\tau}$ we have
$$L_{\s} \otimes M_{\tau} \simeq (L \otimes \s^* M )_{\tau\s}$$
so as left $\oy$-modules we have 
\begin{equation}
L_{\s}^{\otimes j} = L \otimes \s^* L \otimes \s^{2*}L \otimes \ldots \otimes \s^{(j-1)*} L .
\label{eLsj}
\end{equation}
To understand our relation $\phi$, we need to know that a non-zero morphism $\psi:L_{\s}\lm M_{\tau}$ of invertible bimodules exists only if $\s=\tau$ in which case they are given by non-zero morphisms $\tilde{\psi}:L \lm M$ of line bundles. The map $\tilde{\psi}$ gives the map of left modules and the condition $\s = \tau$ ensures that it is compatible with the right module structure too. 
Any line bundle $L$ on $Y$ can be considered as the $\oy$-bimodule $L_{\id}$ so the left and right $\oy$-module structures are the same. Hence our relation $\phi: L_{\s}^{\otimes e} \lm \oy$ is just a non-zero morphism
$$\tilde{\phi}:L \otimes \s^* L \otimes \s^{2*}L \otimes \ldots \otimes \s^{(e-1)*} L \lm \oy .$$
The multiplication map is easy to define too. Given a $G$-invariant open set $U$ and ``left'' sections $s\in L_{\s}^{\otimes i}, t \in L_{\s}^{\otimes j}$ over $U$, we have the product $st$ given by $s \otimes \s^{i *} t$, a left section of $L_{\s}^{\otimes i + j}$ over $U$, or its image in $L_{\s}^{\otimes i + j-e}$ under $\tilde{\phi}$. Let $Z:= Y/G$, the scheme-theoretic quotient which by our assumptions is a quasi-projective variety. Then our cyclic algebra can also be viewed as a finite sheaf of algebras on $Z$.

Suppose now that $e$ is the order of $\s$ and that the relation $\phi:L_{\s}^{\otimes e} \xrightarrow{\sim} \oy$ is an isomorphism. In the commutative case, the latter condition means that the cyclic cover is \'etale, a fact we prove later in the non-commutative case. We now simply write $A$ for the non-commutative cyclic algebra $A(Y; L_{\s}, \phi)$. 

The cyclic algebra $A$ is naturally $G$-graded so there is a natural action of the dual group $G^{\vee}=\langle \s^{\vee}\rangle$ on it. We may thus form the skew polynomial ring 
$$A[u;\s^{\vee}] = A \oplus Au \oplus Au^2 \oplus \ldots $$ 
where for any section $a$ of $A$ we have the skew-commutation relation $ua = \s^{\vee}(a)u$. Note that $Au$ is an $A$-bimodule and in fact, isomorphic to $A_{\s^{\vee}}$ i.e. as a left module it is $A$ and the right module structure is skewed through by $\s^{\vee}$. Recall,
\begin{prop} (\cite[proposition~10.1]{C})
There is a natural action of $G$ on $B:= \send _Y A$ such that $A = B^G$. 
\end{prop}

The action is defined as follows. Let $L_j$ denote the left module structure on $L_{\s}^{\otimes j}$ as in equation~(\ref{eLsj}). Then as sheaves on $Y$ we have $\s^*\shom_Y(L_i,L_j) \simeq \shom(L_{i+1},L_{j+1})$ and summing over $i,j$ gives the action of $\s:\s^*B \xrightarrow{\sim} B$. 

\begin{lemma} \label{lfflat}
\begin{enumerate}
\item $A$ is a faithfully flat left and right $\oy$-module. 
\item $B \simeq A[u;\s^{\vee}]/(u^e-1)$ which is naturally $G^{\vee}$-graded with graded decomposition
$$ B = A \oplus Au  \oplus \ldots \oplus Au^{e-1} .$$
\item If $p:B \lm A$ denotes projection onto the degree 0 component, then $p$ induces an isomorphism of $(B,A)$-bimodules $B \simeq B^*:=\Hom_A(B,A)$. (This is similar to \cite[theorem~2.15~v)]{A}). 
\end{enumerate}
\end{lemma}
\textbf{Proof.}
For part i), we need only note that $A = \oy \oplus L_{\s} \oplus \ldots \oplus L_{\s}^{e-1}$, so is locally free as an $\oy$-module on both the left and right.

We now prove part ii). If we view the $\oy$-module $A$ as a row vector $(\oy \ L_{\s} \ldots L_{\s}^{\otimes e-1})$ then $B$ can be viewed as a matrix algebra 
\[B = \begin{pmatrix}
  \oy & L_{\s} & \hdots & L_{\s}^{\otimes e-1} \\
  L_{\s}^{-1} & \oy & \ddots & \vdots \\
  \vdots & \ddots & \ddots & \vdots \\
  L_{\s}^{\otimes 1-e} & & & \oy
  \end{pmatrix}
\]
Let $\z$ be a primitive $e$-th root of unity and $u$ be the global section of $B$ given by the diagonal matrix
\[u =\begin{pmatrix}
  1 & 0 & \hdots & 0 \\
  0 & \z & \ddots & \vdots \\
  \vdots & \ddots & \ddots & 0 \\
  0 & \hdots & 0 & \z^{e-1}
  \end{pmatrix}\]
The action of $\s$ on $B$ is such that $u$ is a $\z$-eigenvector so there is an eigenspace decomposition
$$B = A \oplus A u \oplus \ldots \oplus A u^{e-1} $$
of $A$-bimodules. This yields ii). The isomorphism in part iii) is given by $b\mapsto p(-b)$.
\qed

\begin{thm} \label{tmorita}
There is a Morita equivalence between $B$ and $\oy$ given by the functor $\Hom_B(A,-)$. If $M$ is an $A$-module, then under this functor, the $B$-module $B \otimes_A M$ corresponds to the underlying $\oy$-module $\ _Y M$. 
\end{thm}
\textbf{Proof.} 
The first assertion holds since $A$ is a vector bundle on $Y$. We use the lemma to see   
\[\Hom_B(A,B \otimes_A M) = \Hom_B(A,B^* \otimes_A M) = \Hom_B(A,\Hom_A(B,M)) = \Hom_A(B \otimes_B A,M) = M\]
as desired.
\qed

\begin{defn}\label{detale}
Let $C \lm C'$ be an homomorphism of algebras (in some abelian category). We say that $C'/C$ is {\em \'etale} if the multiplication map $C' \otimes_C C' \lm C'$ has a $C'$-bimodule splitting. 
\end{defn}

\begin{prop}\label{petale}
The natural morphisms $\oy \hookrightarrow A$ and $A \hookrightarrow B$ are \'etale.
\end{prop}
\textbf{Proof.} 
We show that $\oy \lm A$ is \'etale. Note that 
\[ A \otimes_Y A = \bigoplus_{i,j = 0}^{e-1} L_{\s}^{\otimes i+j} .\]
The $A$-bimodule splitting is easy to construct generically and we do this case first.  Write $K$ for $k(Z)$ so generically, $A$ is $A_K = K[z;\s]/(z^e - \beta)$ for some $\beta \in K^G$ where $K[z;\s]$ denotes the skew polynomial ring (see \cite[example~3.3]{C} for details). Let $\nu_K:A_K \lm A_K \otimes_K A_K$ be the left $A_K$-module module which sends 1 to  
\[c:= \frac{1}{e}(1 \otimes 1 + \sum_{i=1}^{e-1} \beta^{-1}z^i \otimes z^{e-i}) .\]
Note that if $\mu_K:A_K \otimes_K A_K \lm A_K$ is the multiplication map, then $\mu_K(c)=1$. 
To show that $c$ induces a bimodule splitting of $\mu_K$, we need only check that $c$ commutes with multiplication by elements of $A_K$. Now for $\alpha \in K$, it is clear that $\alpha c = c \alpha$. Also, one computes readily that $zc = c z$ so the generic case is settled.

The splitting $\nu_K$ globalises to $\nu:A \lm A \otimes_Y A$ as follows. Given $i,j$ write $n = i+j$ if $i+j<e$ and $n = i+j-e$ otherwise. In the first case, we define $\nu:L_{\s}^{\otimes n} \lm L_{\s}^{\otimes i+j}$ as $\frac{1}{e}$ times the natural isomorphism and in the second case as $\frac{1}{e}$ times $L_{\s}^{\otimes n}  \xrightarrow{\phi^{-1}} L_{\s}^{\otimes n+e} \lm L_{\s}^{\otimes i+j}$. Summing these maps gives a globalisation of $\nu_K$ and thus shows $\oy \lm A$ is \'etale.

By the lemma, $B = A[u;\s^{\vee}]/(u^e-1)$ so the proof that $A \lm B$ is \'etale is similar to the generic case above.
\qed

\section{Line bundles on quaternion orders} \label{slineb}

Assume from now on that $Y$ is a smooth projective surface. Suppose furthermore, that our non-commutative cyclic cover is a double cover $A = \oy \oplus L_{\s}$ which is \'etale in the sense that the relation $\phi$ is an isomorphism. We shall further assume that its ring of fractions $k(A)$ is a division ring which corresponds to the fact that $L$ is non-trivial in $H^1(G, \Pic Y)$ so is not a 1-coboundary. We know from \cite[theorem~3.6]{C} that $A$ is a maximal order in $k(A)$.

For an $A$-module $M$, we define its {\em rank} to be $\rk M := \dim_{k(A)} k(A) \otimes_A M$. We say that $M$ is a {\em vector bundle} if it is locally projective as an $A$-module. It is a {\em line bundle} if its rank is also one. A vector bundle is {\em simple} if its endomorphism ring is $k$. The following is well known.

\begin{lemma} \label{llinebsimple}
Consider an order $B$ in a division ring on a normal projective variety. Let $M$ be a line bundle over $B$. Then $M$ is simple over $B$. 
\end{lemma}
\textbf{Proof.} (We learnt this argument from \cite{AdJ}). $\End_B M$ is a finite dimensional algebra over $k$. It is a domain being a subalgebra of $k(B)$ so must be $k$. \qed

\begin{prop} \label{psplitlineb}
Let $C,C'$ be divisors on $Y$. Then $A \otimes_Y \calo(C) \simeq A \otimes_Y \calo(C')$ if and only if $\calo(C) \simeq \calo(C')$ or $\calo(C) \simeq L_{\s} \otimes_Y \calo(C')$. 
\end{prop}
\textbf{Proof.} This follows from Atiyah's Krull-Schmidt theorem for vector bundles \cite{At}. \qed

We seek to determine all $A$-line bundles. The simplest examples of line bundles are those of the form $A \otimes_Y N$ where $N$ is a line bundle on $Y$. The above proposition shows these are classified by $\Pic Y$ modulo the action $L_{\s} \otimes_Y -$. This is already an interesting object. The next theorem shows there is a dichotomy between line bundles of this form and ones which are simple over $Y$. Recall that the dual group $G^{\vee} = \langle \s^{\vee}\rangle$ acts on $A$ as in lemma~\ref{lfflat}, so we may pullback $A$-modules via $\s^{\vee}$.  

\begin{thm} \label{tdichot}
Let $M$ be a line bundle on $A$. Then either $M$ is simple over $Y$ or it has the form $A \otimes_Y N$ where $N$ is a line bundle on $Y$. The latter occurs precisely when $M \simeq \s^{\vee *} M$. 
\end{thm}
\textbf{Proof.} Recall lemma~\ref{lfflat} shows that $B:= \send_Y A = A \oplus A_{\s^{\vee}}$. Using the Morita equivalence between $B$ and $\oy$ as made explicit in theorem~\ref{tmorita} we see 
\begin{align*}
\End_Y M & = \Hom_B (B \otimes_A M,B \otimes_A M)  \\
& = \Hom_A(M,M \oplus \s^{\vee *} M)  \\
& = k \oplus \Hom(M,\s^{\vee *} M) 
\end{align*}
by simplicity of $M$ over $A$. If $\Hom(M,\s^{\vee *} M) = 0$ then $M$ is simple over $Y$ so suppose it is non-zero. Now $\s^{\vee *}: A -\mbox{Mod} \lm A-\mbox{Mod}$ is a category equivalence so 
\[ \Hom(M,\s^{\vee *} M) = \Hom_A(\s^{\vee *}M,M).\]
But $M,\s^{\vee *} M$ are simple bundles over $A$ which are critical so this is non-zero if and only if $M \simeq \s^{\vee *}M$. Thus $\End_Y M \simeq k \times k$ as vector spaces. To see this also an isomorphism of rings, note that the above computation shows there is a non-zero $B$-module map $\tau:B \otimes_A M \lm B \otimes_A M$ which swaps $M$ and $\s^{\vee *} M$. Hence there are primitive orthogonal idempotents in $\End_Y M$ corresponding to $\frac{1}{2}(1 \pm \tau)$. These idempotents yield a direct sum decomposition $M = N \oplus N'$ where $N,N'$ are line bundles on $Y$. We consider the multiplication map $m: L_{\s} \otimes_Y M  \lm M$ and the composite map 
\[ \oy \otimes_Y M \simeq L_{\s} \otimes_Y L_{\s} \otimes_Y M \xrightarrow{1 \otimes m} L_{\s} \otimes_Y M \xrightarrow{m} M \]
which by associativity is an isomorphism. This shows $m$ is surjective and it must be injective as the kernel must have rank 0. We thus have the following isomorphism of $\oy$-modules $L_{\s} \otimes_Y N \oplus L_{\s} \otimes_Y N' \simeq N \oplus N'$. Atiyah's Krull-Schmidt theorem shows that $L_{\s} \otimes N$ is isomorphic to $N'$ or $N$. In the first case we see $M \simeq A \otimes_Y N$. In the second case we find $L \simeq N \otimes_Y \s^* N^{-1}$ which contradicts the fact that $L$ is not a 1-coboundary. 
\qed

\vspace{2mm}

One approach to constructing bundles on $A$ is to start with even rank bundles on $Y$ and impose $A$-module structures on it. For rank two bundles on $Y$, we see the possibilities are fairly limited.

\begin{prop} \label{pAstr}
Let $M$ be a rank two vector bundle on $Y$. If $M$ is split, then there is at most one $A$-module structure on $M$ up to isomorphism. If $M$ is simple (over $Y$), then there are two possibilities. Either $L_{\s} \otimes_Y M$ is not isomorphic to $M$ in which case $M$ cannot be given the structure of an $A$-module. If $L_{\s} \otimes_Y M \simeq M$ then there are exactly two possible $A$-module structures on $M$.
\end{prop}
\textbf{Proof.} We do the simple case only since the other is clear. If $M$ is an $A$-module then the multiplication map $m: L_{\s} \otimes_Y M \lm M$ is an isomorphism by the argument in the proof of theorem~\ref{tdichot}. Suppose conversely that we are given a simple rank two bundle $M$ with an $\oy$-module isomorphism $m: L_{\s} \otimes_Y M \lm M$. All other isomorphisms are given by changing $m$ by a multiplicative constant $\alpha\neq 0$. Note $m$ defines an $A$-module structure if and only if it satisfies the cocycle condition, that is, 
\[ M \simeq L_{\s} \otimes_Y L_{\s} \otimes_Y M \xrightarrow{1 \otimes m} L_{\s} \otimes_Y M \xrightarrow{m} M \]
is the identity. Changing $m$ by $\alpha$ alters this composite by $\alpha^2$ so there are precisely two values of $\alpha$ which gives valid descent data. 
\qed

\vspace{2mm}

For a line bundle $M$ over $A$, we consider the question of semistability of the underlying $\oy$-module. We fix some ample $G$-invariant divisor $H$ on $Y$. Let $V$ be a vector bundle on $Y$. We define its slope (with respect to $H$) to be 
\[\mu(V):= \frac{c_1(V).H}{\rk V} \]
We say that $V$ is {\em $H$-semistable} if for any subbundle $U < V$ we have $\mu(U) \leq \mu(V)$. 
\begin{lemma}  \label{lslopeL}
If $L=\oy(E)$ then $E.H=0$.
\end{lemma}
\textbf{Proof.}
Since $L$ is a 1-cocycle, we know $E \sim - \s^* E$. Hence $E.H = -\s^* E.H = - E. \s^* H = - E.H$ so $E.H = 0$. 
\qed
\begin{prop}  \label{pmuss}
If $M$ is a line bundle over $A$ then the underlying $\oy$-module $_Y M$ is $H$-semistable. 
\end{prop}
\textbf{Proof.} Let $N<M$ be a line bundle over $Y$. Since any $A$-submodule of $M$ must be rank 2 over $Y$, we see that $L_{\s} \otimes_Y N \cap N = 0$. Now the lemma ensures
\[ c_1(L_{\s} \otimes_Y N).H = (c_1(L) + \s^* c_1(N)).H = c_1(L).H + c_1(N).H = c_1(N).H \]
Hence
\[ \mu(N) = \mu(N \oplus L_{\s} \otimes_Y N) \leq \mu(M) .\] \qed

\vspace{2mm}

In the next section we show we don't have (Gieseker) semistability in general. 

\section{Some exotic del Pezzo orders} \label{sdPcase}

In \cite{CK}, non-commutative analogues of del Pezzo surfaces, dubbed del Pezzo orders were classified. The generic example was a maximal order on $\PP^2$ ramified on a cubic. There were however some exotic families such as the quaternion order on $\PP^2$ ramified on a smooth quartic. An explicit construction of these (up to Morita equivalence) was given via non-commutative cyclic covers in \cite[section~6]{C}. In this section, we review this example and show that its line bundles are not in general (Gieseker) semistable. The following sections will study the moduli of line bundles for these orders. 

For the rest of this paper, we let $\pi:Y\lm \PP^2$ be the double cover of $\PP^2$ ramified on some smooth quartic. Thus $Y$ is isomorphic to the blowup of $\PP^2$ at 7 points in general position. It has 56 exceptional curves, 2 lying over each of the 28 bitangents of the plane quartic. Let $E,E'$ be two disjoint exceptional curves on $Y$ so that $L:=\oy(E - E')$ is a 1-cocycle and hence $L_{\s}^{\otimes 2} = L \otimes \s^* L \simeq \oy$. If we let $\phi: L_{\s}^{\otimes 2} \lm \oy$ be any isomorphism, then it satisfies the overlap condition and $A = \oy \oplus L_{\s}$ is a cyclic algebra. Changing the relation $\phi$ yields an isomorphic order. The order is maximal and ramified on the same quartic that $Y$ is ramified on. 

We can carry out the above constructions for the case when $E,E'$ are exceptional curves which intersect. The calculations will be different but the cyclic algebra $A(Y;\oy(E-E')_{\s})$ is Morita equivalent to
$$ \oy(E') \otimes_Y A(Y;\oy(E-E')_{\s}) \otimes_Y \oy(-E') \simeq  A(Y;\oy(E-\s E')_{\s})  $$
and, since $E\cap\s E'=\varnothing$, the line bundles can be computed from the disjoint case. It is for this reason that we confine ourselves to the case where $E,E'$ are disjoint.

One nice feature of del Pezzo surfaces is that for the moduli of stable bundles, the obstruction to smoothness vanishes. We wish to recall in what sense our cyclic algebra $A$ is also del Pezzo, and hence, exhibits similar behaviour. First, as suggested by the adjunction formula, we define the {\em canonical bimodule} to be $\w_A:=\shom_Z(A,\w_Z)$ which is an $A$-bimodule. One readily obtains Bondal-Kapranov-Serre duality, that is, for $A$-modules $M,N$, there are natural isomorphisms
$$ \Ext^i_A(M,N) \simeq \Ext^{2-i}_A(N,\w_A \otimes_A M)^* .$$

Fortunately, $\w_A$ is completely computable in our case. Let $R \subset Y$ denote the ramification curve of $\pi:Y \lm \PP^2$ which is the inverse image of our smooth quartic. Note that since $R$ is $G$-invariant, $A \otimes_Y \oy(R)$ is naturally an $A$-bimodule. Let $H \in \Div Y$ be the inverse image of a general line in $\PP^2$. Then we have
\begin{prop} \label{pwA}
The canonical bimodule $\w_A= A \otimes_Y \oy(-H)$. In particular, the coarse moduli scheme of $A$-line bundles is smooth.
\end{prop}
\textbf{Proof.}
Putting together \cite[proposition~5(1)]{CK} and \cite[proof of proposition~4.5]{C}, we find that
$$ \w_A \simeq \w_Z \otimes_Z A \otimes_Y \oy(R) \simeq A \otimes_Y \oy(-3H+R) \simeq A \otimes_Y 
\oy(-H)  .$$
We proceed to compute the obstruction to smoothness \cite[lemma~3.2]{HoSt}. Let $M$ be an $A$-line bundle. Then
$$ \Ext^2_A(M,M) \simeq \Hom_A(M, \w_A \otimes_A M)^* \simeq \Hom_A(M, \oy(-H) \otimes_Y M)^* = 0 .$$
The proposition is proved.
\qed

Note that $H$ is ample so in some sense $\w_A^{-1}$ is also ample and $A$ is del Pezzo. The interested reader may consult \cite[definition~7]{CK} for the actual definition of a del Pezzo order. 
To finish off this section we observe the following

\begin{prop}  \label{pAnotss}
The $\oy$-module $A$ is not semistable.
\end{prop}
\textbf{Proof.}
It suffices to show that $\chi(\oy) \neq \chi(\oy(E-E'))$ so their Hilbert polynomials differ too. 
Now the canonical divisor $K$ is $G$-invariant so 
\[2(\chi(\oy(E-E')) - \chi(\oy)) = (E-E').(E-E'-K) = E^2 + E'^2 = -2 \]
by Riemann-Roch.
\qed

\section{Possible Chern classes}   \label{schern}

We continue our notation from the previous section and seek to determine all the possible Chern classes of $A$-line bundles on our del Pezzo order. Recall $Y$ is the blowing up of $\PP^2$ at seven points in general position. By contracting exceptional curves $E_1,\ldots, E_7\subset Y$ to $p_1,\ldots,p_7\in \PP^2$ appropriately, we may assume that $E = E_1$ and $\s E'$ is the strict transform of the line through $p_1,p_2$.

The possible first Chern classes of $A$-line bundles are rather limited. To describe them, let $H \in \Div Y$ be the inverse image under $\pi:Y \lm \PP^2$ of a general line. 
\begin{prop}  \label{ppossibleChern}
If $M$ is an $A$-line bundle then there is some $n \in \Z$ such that $c_1(M) = L_{\s} \otimes_Y \oy(nH)$. 
\end{prop}
\textbf{Proof.} 
We know that $L_{\s} \otimes_Y M \simeq M$ so taking first Chern classes of both sides shows $L^{\otimes 2} \simeq \det M \otimes_Y (\s^* \det M)^*$. To solve this equation for $\det M=c_1(M)$, note first that $c_1(M)=L$ is a solution since $L$ is a 1-cocycle. All other solutions differ by an element of $(\Pic Y)^G$ so let $\oy(D) \in (\Pic Y)^G$. Then
$$ 2D \sim D + \s^* D \sim \pi^*\pi_* D \in \pi^* \Div \PP^2 .$$ 
But $\pi^*\Pic \PP^2$ is a primitive lattice in $\Pic Y$ so $D \in \pi^*\Pic \PP^2$ and $(\Pic Y)^G = \pi^* \Pic \PP^2$.
\qed \vspace{2mm}

It is easy to see that all such Chern classes actually arise. We know from proposition~\ref{pmuss} that any $A$-line bundle is $H$-semistable so we may apply Bogomolov's inequality to bound the second Chern class. Namely, if $\Delta = 4c_2 - c_1^2$ (which is a scalar multiple of the usual discriminant) then $\Delta \geq 0$. The optimal bound is below.

\begin{prop}  \label{pBogup}
Any $A$-line bundle satisfies $\Delta > 0$.
\end{prop}
\textbf{Proof.} Let $M$ be an $A$-line bundle. Since $\Delta$ is unaffected by tensoring with a line bundle, we may assume that $c_1(M) = L_{\s} \otimes_Y \oy(nH)$ where $n=0$ or 1. If $n=0$ then $c_1^2 = -2$ so integrality of Chern classes ensures $\Delta > 0$. We shall thus assume $n=1$ and $\Delta = 0$ which amounts to $c_2(M) = 0$. Let $M_1 := \oy(-E) \otimes_Y M$. Note that 
\[  c_1(M_1) = c_1(M) - 2E, c_2(M_1) = c_2(M) - c_1(M).E + E^2 = -1  .\]
Then Riemann-Roch on $Y$ gives
\[ \chi(M_1) = 2 + \frac{1}{2}(\s E' - E).(\s E' - E + H) + 1 = 1.  \]
Now $M_1^*\otimes_Y \oy(-H)$ is $H$-semistable with first Chern class $E-\s E' - 2H \sim -\s E - \s E' - H$. Hence 
\[ h^2(M_1) = h^0(M_1^* \otimes_Y \oy(-H)) = 0  \]
so $H^0(M_1) \neq 0$ (we can also conclude $h^2 = 0$ from the fact that $M_1^* \otimes_Y \oy(-H)$ has negative slope). There is consequently, an embedding $\oy(E) \hookrightarrow M$ of coherent sheaves on $Y$ which in turn gives an embedding of $A$-modules $\phi: A \otimes_Y \oy(E) \hookrightarrow M$. But $c_1(A \otimes_Y \oy(E)) = c_1(M)$ so $\phi$ must be an isomorphism. This contradicts the fact that $A \otimes_Y \oy(E)$ and $M$ have different second Chern class.

\begin{thm} \label{tchern2} 
Let $n \in \Z$ and $m \in \Z$ be so that 
$$\Delta(n,m) := 4m - c_1(L_{\s} \otimes_Y \oy(nH))^2 > 0.$$ 
Then there is a $Y$-split $A$-line bundle $M = A \otimes_Y \oy(D)$ with 
\[ c_1(M) = L_{\s} \otimes_Y \oy(nH), c_2(M) = m  .\]
\end{thm}
\textbf{Proof.} 
Let $f:Y \lm \PP^2$ be the blowing down which contracts $E_1,\ldots, E_7$. We write $F$ for the inverse image in $Y$ of a generic line in $\PP^2$ so that $H \sim 3F - \sum E_i$. Note that $E_i + \s E_i \sim H$ and 
$$ 2H \sim H + \s H \sim (3F - \sum E_i) + (3 \s F - \sum (H-E_i)) \sim 3(F + \s F) - 7H$$
so $F + \s F \sim 3H$. If 
\[ D = n_0 F + \sum_{i=1}^7 n_i E_i  \]
then for $M = A \otimes_Y \oy(D)$ we compute
\[ c_1(M) = D + E - E' + \s D = E - E' + (3n_0 + \sum_{i>0} n_i)H  .\]
This gives 
\begin{equation} 
n =  3n_0 + \sum_{i>0} n_i .
\label{en}
\end{equation}
Tensoring by $\oy(H)$ preserves the discriminant $\Delta$ so we may restrict to the cases $n=0,1$. 

Note that $2n^2 = (nH)^2 = (D+ \s D)^2$
which gives the useful formula 
\[ D. \s D = n^2 - D^2   .\]
Also, $E'$ is the strict transform of the unique conic through $p_3,\ldots, p_7$ so 
\[ c_2(M) = D.(E - E' + \s D) = D.(E_1 -2F + E_3 +E_4 + E_5+E_6 + E_7) + n^2 - D^2 .\]
The possible values this can take is determined by a simple number theoretic problem. We will use 
\begin{lemma} \label{lsquares} 
Any odd positive integer can be written as a sum of squares $\sum_{j=1}^4 t_j^2$ for integers $t_j$ satisfying $\sum t_j = 1$. 
\end{lemma} 
\textbf{Proof.} This is part of \cite[exercise~15, chapter~9]{Rose}. 

We now complete the proof of the proposition by dividing into cases.

\underline{$n=1$:} For any $m > 0$, we need to solve (\ref{en}) and 
\[ m = 1-n_0^2 -2n_0+ (n_1^2 -n_1) + n_2^2 + \sum_{i \geq 3} (n_i^2 - n_i).\]
For $m=1$ we obtain a solution on setting $n_0=1,n_1=-1,n_2=-1,n_3 = \ldots =n_7 = 0$ whilst if $m=2$ we have the solution $n_0 = -2,n_1= \ldots = n_7 = 1$. For $m>2$, 
we seek solutions where $n_0 = -2$ and $n_4,n_5,n_6,n_7$ vary subject to $n_4 + \ldots + n_7 = 5$ so that  
\[ c_2(M) = 2 + n_2^2 + (n_1^2 -n_1) + (n_3^2 -n_3) +\sum_{i \geq 4} (n_i - 1)^2.\] 
If also $n_1 = n_3=1,n_2=0$ then (\ref{en}) holds and we see on applying the lemma that $c_2(M)$ can be any odd integer $m\geq 3$. If instead $n_1 = n_2 = 1, n_3 = 0$ we see that (\ref{en}) holds and $m$ can be any even integer $\geq 4$. 

\underline{$n=0$:} For any $m \geq 0$, we need to solve (\ref{en}) and 
\[ m = -n_0^2 - 2n_0 + n_2^2 + (n_1^2 -n_1)  + \sum_{i \geq 3} (n_i^2 - n_i).\]
Note the $m=0$ case is solved with all $n_i=0$. For $m=1$ we get the solution $n_2= -1,n_3=1$ and all other $n_i=0$. For $m >1$, we let $n_0 = -2$ and $n_4,n_5,n_6,n_7$ vary subject to $n_4 + \ldots + n_7 = 5$ so that  
\[ c_2(M) = 1 + n_2^2 + (n_1^2 -n_1) + (n_3^2 -n_3) +\sum_{i \geq 4} (n_i - 1)^2.\]
If $m$ is odd, we obtain a solution on setting $n_2 = 1,n_1 = n_3 = 0$ while if $m$ is even we obtain a solution with $n_1 = 1,n_2 = n_3 = 0$.

The theorem is now proved. 
\qed  \vspace{2mm}

\section{Moduli Spaces in the Case of Minimal Second Chern Class}

We continue the notation of the last two sections. Our aim in this section is to compute moduli spaces of $A$-line bundles when the second Chern class is minimal. Fortunately, thanks to the work of Simpson, Hoffmann-Stuhler and Lieblich (see \cite{Simp},\cite[section~2]{HoSt},\cite{Lieb}), we know that the commutative theory of coarse moduli schemes for torsion-free sheaves pushes through to our case. In particular, since we are only concerned with $A$-line bundles, semsistability considerations do not arise. 

As in proposition~\ref{pBogup}, we may assume that $c_1 = L \otimes_Y \oy(nH)$ where $n = 0$ or 1. The corresponding minimal second Chern classes are, according to proposition~\ref{pBogup}, 0 and 1 respectively. We look at the $n=0$ case first.

\begin{prop}  \label{pnis0minc2}
Let $M$ be an $A$-line bundle with $c_1 = L$ and $c_2 = 0$. Then $M \simeq A$. In particular, the coarse moduli scheme of such $A$-line bundles is a point. 
\end{prop}
\textbf{Proof.} If $H^0(M) \neq 0$ then $A$ embeds in $M$ so a comparison of first Chern classes shows that $A \simeq M$ as desired. Note $A$ and $M$ have the same Chern classes so Riemann-Roch tells us that $\chi(M) = \chi(A) = 1$. It thus suffices to show that $H^2(M) = 0$. Assume this is not the case. Recall from proposition~\ref{pwA} that $\w_A = A \otimes_Y \oy(-H)$ so $c_1(\w_A) = E-E' -2H$. Hence $c_1(\w_A) - c_1(M)$ cannot be effective. This means that there is no embedding $M \hookrightarrow \w_A$. 
Serre duality shows that $H^2(M) = \Hom_A(M,\w_A) = 0$ so taking into account the Euler characteristic we see that $H^0(M) \neq 0$ as desired. 

For the last statement, we compute the tangent space to the moduli space with the aid of Leray-Serre and \cite[lemma~3.1]{HoSt},
\[\Ext^1_A(A,A) = \Ext^1_Y(\oy, \oy \oplus \oy(E-E')) = H^1(\oy \oplus \oy(E-E')) = 0 .\]
Hence $A$ is rigid and the moduli scheme is a point. 

\qed \vspace{2mm}

The other minimal Chern class case occurs when $c_1 = \s E' + E$ and $c_2 = 1$. This includes the $A$-line bundle $A \otimes_Y \oy(\s E') = \oy(\s E') \oplus \oy(E)$. We seek to construct a 1-parameter family of $A$-line bundles which has this line bundle as a member. In fact it will have 6 members which are $Y$-split.

To motivate the construction, we describe first the $A$-line bundles which form the individual members of the family. First note that $(E + \s E')^2 = 0$ and standard cohomology arguments now show that $|E + \s E'|$ is a basepoint free pencil on $Y$. Pick a (closed) point $p \in Y$ which belongs to a unique curve $F(p)$ in $|E + \s E'|$. 

To understand this linear system better, recall that we can contract exceptional curves $E_1,\ldots, E_7$ to $p_1,\ldots,p_7$ so that $E = E_1$ and $\s E'$ is the strict transform of the line through $p_1,p_2$. The generic member of $|E + \s E'|$ is the strict transform of a line through $p_2$. There are 6 other members which are pairs of intersecting exceptional curves. 

Let $H \in \Div Y$ be the inverse image under $\pi:Y \lm \PP^2$ of a general line. From \S~\ref{schern}, we know that a possible first Chern class for an $A$-line bundle is $L_{\s} \otimes_Y \oy(H) \simeq \oy(E + \s E')$. We will construct some line bundles with this as its first Chern class. If $I_p \triangleleft \oy$ denotes the ideal sheaf at $p$ and $k(p)$ the corresponding skyscraper sheaf, then the long exact sequence in cohomology arising from the exact sequence 
\[ 0 \lm I_p\oy(E + \s E') \lm \oy(E + \s E') \lm k(p) \lm 0 \]
shows $\Ext^1_Y(I_p\oy(E + \s E'),\oy) = \Ext^2_Y(k(p),\oy) = k$. Consequently, by Serre's lemma (see \cite[lemma~5.1.2]{OSS} or \cite[chapter~5,section~4, lemma page~724]{GH}), there is a unique rank two vector bundle $M = M(p)$ on $Y$ which fits in a non-split exact sequence
\begin{equation}\label{eMp} 
0 \lm \oy \lm M(p) \lm I_p\oy(E + \s E') \lm 0   
\end{equation} 
The sequence shows that $H^0(M(p)) = k^2$. 
To determine which $M(p)$ are isomorphic, we compute the possible zeros of non-zero sections. 

\begin{lemma} \label{lzeros}
Let $F \in |E + \s E'|$ be the curve in the linear system which contains $p$. If $F$ is irreducible then the zeros of (non-zero) global sections of $M(p)$ is a point which varies over the whole of $F$. If $F = F_1 + F_2$ where $F_1,F_2$ are exceptional curves, then we have an exact sequence 
\[0 \lm \oy(F_i) \lm M(p) \lm \oy(F_j) \lm 0 \]
where $p \in F_j$ and $F_i$ is the other curve. Moreover, the sequence is non-split unless $p = F_1 \cap F_2$ in which case it does split. 
\end{lemma}
\textbf{Proof.}
Suppose first that $F$ is irreducible. Let $s:\oy \lm M$ represents a non-zero section. If it is the map in the exact sequence (\ref{eMp}) above, its zero set is $p$. Otherwise, we can compose to obtain an injection $\oy \lm I_p\oy(E + \s E') \simeq I_p\oy(F)$. Since $F$ is irreducible, the image of this map must be $\oy(F-F')$ for some effective $F'\sim F$  such that $p \in F'$. But $|F|$ is a basepoint free pencil so $F = F'$. Hence the zeros of $s$ must lie on $F$ and looking at the closed fibre at $p$ we see it is not $p$. Varying the global section $s$ we see its zero set is a point which varies over $F$. 

Suppose now that $F=F_1+ F_2$ is reducible as above with say $p \in F_2$. Then $\Hom_Y(\oy(F_1),I_p\oy(F)) \neq 0$ so from (\ref{eMp}), we see there exists an embedding of $\oy(F_1)$ into $M(p)$. That $M(p)$ is an extension of line bundles as above now follows from Chern class computations. \qed

We check semistability of the $M(p)$.

\begin{prop} \label{pMpsemis}
If $F(p)$ is reducible, then $M(p)$ is (Gieseker) semistable and hence $H$-semistable. If $F(p)$ is irreducible then $M(p)$ is $H$-stable and hence stable.
\end{prop}
\textbf{Proof.}
First note that Riemann-Roch gives the Hilbert polynomial for $\oy(D)$ as 
\[ \chi(\oy(D +nH)) = \frac{1}{2}n^2 H^2 + \frac{1}{2}n(2D.H - K.H) + \chi(\oy) + \frac{1}{2}(D^2 - D.K) .\] 
If $F(p) = F_1 + F_2$ then without loss of generality we may suppose $M(p)$ sits in an exact sequence 
\[0 \lm \oy(F_1) \lm M(p) \lm \oy(F_2) \lm 0 \]
Semistability follows now from the fact that $\oy(F_1),\oy(F_2)$ have identical Hilbert polynomial. 

Suppose now that $F = F(p)$ is irreducible. By moving $p$ in $F$, we may assume that $p$ does not lie on an exceptional curve. The slope of $M=M(p)$ is 
\[ \mu(M) = \frac{1}{2}c_1(M).H = \frac{1}{2}(E+\s E').H = 1 .\]
Let $N$ be a line bundle which embeds in $M$. Suppose it has positive slope so that it must embed in $I_p\oy(F)$. Then $N \simeq \oy(F-D)$ where $D$ is an effective divisor containing $p$. Since $D$ is not exceptional, $\pi_* D$ is not a line (with multiplicity one) so $D.H \geq 2$. Hence
\[\mu(N) = (F-D).H \leq 0 < \mu(M) \]
proving $H$-stability. 
\qed \vspace{2mm} 

The next result gives a simple criterion for a rank two bundle on $Y$ to be of the form $M(p)$. 
\begin{schol}  \label{scMptype}
Any $H$-semistable rank two vector bundle $V$ on $Y$ with $\dim_k H^0(Y,V) \geq 2$ and Chern classes $c_1(V) = \oy(E + \s E'),c_2(V) = 1$ is isomorphic to $M(p)$ for some $p$. 
\end{schol}
\textbf{Proof.}
Suppose a non-zero global section of $V$ has isolated zeros. Then it gives rise to an extension as in (\ref{eMp}). Otherwise, the zeros yield an effective divisor $D_1$ such that there is an exact sequence
\[ 0 \lm \oy(D_1) \lm V \lm I\oy(D_2) \lm 0 \] 
where $I$ is some ideal sheaf and $D_2 \in \Div Y$. $H$-semistability means that $D_1.H \leq 1$. The only possibility is that $D_1$ is an exceptional curve. The exact sequence above now shows that $H^0(Y,\oy(D_2)) \neq 0$ so we may as well assume that $D_2$ is effective. Comparing first Chern classes we see that $D_1+D_2 \sim E+\s E'$. Looking at the second Chern class shows that $I = \oy$ so $V \simeq M(p)$ where $p \in D_2$ by lemma~\ref{lzeros}. 
\qed \vspace{2mm} 

We now check when $L_{\s} \otimes_Y M(p) \simeq M(p)$ and so, by proposition~\ref{pAstr}, gives rise to an $A$-line bundle. The first step is to use the previous criterion to show that $L_{\s} \otimes_Y -$ sends bundles of the form $M(p)$ to bundles of the same form. The Chern classes of $L_{\s}\otimes_Y M(p)$ are easy enough to determine from the exact sequence 
\[ 0 \lm \oy(E-E') \lm L_{\s} \otimes_Y M \lm I_{\s(p)}\oy(H) \lm 0 .\]
We see indeed that $c_1(L_{\s} \otimes_Y M(p)) = E + \s E', c_2(L_{\s} \otimes_Y M(p)) = 1$. It also shows that $H^0(Y,L_{\s} \otimes_Y M(p)) = 2$. Finally we have
\begin{prop}\label{pstillsemis}
Let $V$ be an $H$-semistable (resp. $H$-stable) bundle on $Y$. Then $L_{\s} \otimes_Y V$ is also $H$-semistable (resp. $H$-stable). 
\end{prop}
\textbf{Proof.}
If $W < L_{\s} \otimes_Y V$ is a rank $r$ subbundle then $L_{\s} \otimes_Y W$ is a subbundle of $V$ and 
\[c_1(L_{\s} \otimes_Y W).H = (rE-rE' + \s^* c_1(W)).H = c_1(W).H  .\]
The proposition now follows.
\qed \vspace{2mm}

\begin{cor} \label{cfamiliyinv}
For any $p \in Y$ there is some point $q \in Y$ with $M(q) \simeq L_{\s} \otimes_Y M(p)$. 
\end{cor}

Determining which bundles $M(p)$ satisfy $L_{\s} \otimes_Y M(p) \simeq M(p)$ should just be a simple matter of computing zeros of sections and applying lemma~\ref{lzeros}. We were unable to make this approach work so we proceed somewhat more subtly. We check solutions in the $F(p)$ reducible case first. 
This occurs whenever $M(p)$ is $Y$-split. Then we build a $\PP^1$-family of bundles of the form $M(p)$. This gives a rational curve in the coarse moduli scheme of semistable rank two bundles on $Y$. Tensoring by $L_{\s}$ maps this rational curve to itself and also fixes the 6 points corresponding to the reducible members of the linear system $|E + \s E'|$. It thus must fix the whole curve. 

We deal with the $F(p)$ reducible case first. 
\begin{lemma} \label{lredflips}
Let $F_1,F_2$ be exceptional curves such that $F_1+F_2 \sim E + \s E'$. Consider an extension 
\[0 \lm \oy(F_1) \lm M \lm \oy(F_2) \lm 0  .\]
Tensoring by $L_{\s}$ gives the exact sequence
\[0 \lm \oy(F_2) \lm L_{\s}\otimes_Y M \lm \oy(F_1) \lm 0  .\]
In particular, $L_{\s} \otimes_Y M \simeq M$ if and only if $M$ is split in which case it corresponds to the $A$-line bundle $A \otimes_Y \oy(F_1) \simeq A \otimes_Y \oy(F_2)$. 
\end{lemma}
\textbf{Proof.}
Observe
$$ E - E' + \s F_1 \sim E - E' + \s E + E' - \s F_2 \sim H - \s F_2 \sim F_2$$
so $L_{\s} \otimes_Y \oy(F_1) \simeq \oy(F_2)$ and similarly $L_{\s} \otimes_Y \oy(F_2) \simeq \oy(F_1)$. The lemma follows. 
\qed \vspace{2mm}

Note that the non-split extension above is not fixed by $L_{\s} \otimes_Y-$, however, its associated graded is, and hence, so is its image in the coarse moduli scheme of semistable rank two bundles on $Y$. 

We now construct a $\PP^1$-family of $M(p)$'s. Recall that we may assume that $Y$ is the blow-up of $\PP^2$ at the points $p_1,\ldots, p_7$, that $E$ is the exceptional curve above $p_1$ and that $\s E'$ is the strict transform of the line $p_1p_2$. Choose a smooth rational curve $C \subset Y$ with $C.(E+ \s E') = 1, C^2 = 0$. For example, we can take $C$ to be the strict transform of a generic line through $p_3$. We will abuse notation and also let $C$ denote its image under the diagonal embedding $C \lm Y \times C$ and let $I_C$ denote the corresponding ideal sheaf. Since $C$ is smooth rational, $\calo_C(1)$ has its obvious meaning. Note that $C$ is a section to the fibration $|E + \s E'|: Y \lm \PP^1$.

Let $\delta: Y \times C \lm C, \eps:Y \times C \lm Y$ denote the projection maps. For sheaves $M_1$ on $Y$, $M_2$ on $C$ we let $M_1 \boxtimes M_2 := \eps^* M_1 \otimes_{Y \times C} \delta^* M_2$. 

\begin{lemma} \label{lP1family}
There exists a vector bundle $M$ on $Y \times C$ which fits in an exact sequence 
\begin{equation} 
0 \lm \calo_{Y \times C} \lm M \lm I_C \oy(E + \s E') \boxtimes \calo_C(1) \lm 0  .
\label{efamily}
\end{equation}
For any $p \in C$, restricting the above sequence to $Y \times p$ yields the exact sequence 
\[0 \lm \oy \lm M(p) \lm I_p\oy(E + \s E') \lm 0   .\]
In particular, we obtain an injective map of $C$ into the coarse moduli scheme of semistable rank two vector bundles on $Y$ with $c_1 = E + \s E', c_2 = 1$. 
\end{lemma}
\textbf{Proof.} 
We compute the appropriate Ext group using the local-global spectral sequence. To simplify notation in this proof, we shall omit the default subscript $Y \times C$ and only retain the subscript $Y,C$ or otherwise as the case may be. For example, $\w$ will mean the canonical sheaf $\w_{Y \times C}$. 
Note 
\[ \Rshom(I_C \oy(E + \s E') \boxtimes \calo_C(1),\calo) = \Rshom(I_C,\w) \otimes \w^* \otimes 
(\oy(-E - \s E') \boxtimes \calo_C(-1))  \]
Consider the exact triangle 
\[ \Rshom(\calo_C,\w) \lm \Rshom(\calo,\w) \lm \Rshom(I_C,\w)  .\]
Grothendieck duality theory gives $\Rshom(\calo_C,\w) = \w_C [2]$, the $[2]$ denoting shift in cohomological degree. Hence 
\[\shom(I_C,\w) = \w, \sext^1(I_C,\w) = \w_C, \sext^2(I_C,\w) = 0  .\]
Hence 
\begin{align*}
\shom(I_C \oy(E + \s E') \boxtimes \calo_C(1),\calo) & =  \oy(-E - \s E') \boxtimes \calo_C(-1) \\
\sext^1(I_C \oy(E + \s E') \boxtimes \calo_C(1),\calo) & = \w_C \otimes \w^* \otimes 
(\oy(-E - \s E') \boxtimes \calo_C(-1))
\end{align*}
To simplify notation, we abbreviate the two terms above as $\shom,\sext$. Note that $\delta$ restricted to $C$ is an isomorphism so 
\[ \w_C \otimes \w^* = \w_C \otimes \delta^*\w^*_C \otimes \eps^*\w_Y^*  = \eps^*\w_Y^*  \]
Also, $\eps$ restricted to $C$ is an isomorphism so 
\[ \eps^*\oy(-E - \s E')|_C = \calo_C(-E - \s E') = \calo_C(-1) \]
by our hypothesis on $C$. Similarly the hypothesis $C^2 = 0$ gives $K.C=-2$ and hence
\[ \eps^* \w^*_Y|_C = \calo_C(2)   .\]
We may thus simplify 
\[ \sext = \calo_C(2) \otimes_C \calo_C(-1) \otimes_C \calo_C(-1) \simeq \calo_C  .\]
We find that $H^0(Y\times C, \sext) = H^0(C,\calo_C) = k$ and the K\"unneth formula gives $H^2(Y\times C,\shom) = 0$. From the local-global spectral sequence, we see that the non-zero section of $H^0(Y\times C, \sext)$ gives a non-trivial extension class in $\Ext^1(I_C \oy(E + \s E') \boxtimes \calo_C(1),\calo)$. The section clearly generates $\sext$ locally at every point so by Serre's lemma, we get a vector bundle $M$ which fits in the exact sequence (\ref{efamily}) above. Finally, $I_C \oy(E + \s E') \boxtimes \calo_C(1)$ is flat over $C$ so restricting $M$ to any fibre of $\delta: Y \times C \lm C$ yields $M(p)$ as desired.
\qed \vspace{2mm} 

\begin{thm} \label{tsomenonsplit}
If $F(p)$ is irreducible then $L_{\s} \otimes_Y M(p) \simeq M(p)$. In particular, $M(p)$ carries two possible structures as $A$-line bundles. 
\end{thm}
\textbf{Proof.} 
The family in the previous lemma, and that obtained by tensoring by $L_{\s}$ give two maps of $\PP^1$ into the coarse moduli scheme of semistable rank two bundles on $Y$. They are both injective and agree on 6 points so must be the same map. 
\qed \vspace{2mm}

\begin{prop}  \label{pnis1minc2}
Let $M$ be an $A$-line bundle with Chern classes $c_1 = L \otimes_Y \oy(H), c_2 = 1$. Then as $\oy$-modules we have $M \simeq M(p)$ for some $p$. 
\end{prop}
\textbf{Proof.} We can assume that $M$ is $Y$-simple. Proposition~\ref{pmuss} shows that $M$ is $H$-semistable so by scholium~\ref{scMptype}, it suffices to show that $\dim_k H^0(Y,M) \geq 2$. Note 
\[ \chi(M) = \chi(A \otimes_Y \oy(\s E')) = 2  .\]
But $H^2(M) = \Hom_A(M,\w_A)^*$. Comparing first Chern classes of $M$ and $\w_A$ as in the proof of proposition~\ref{pnis0minc2} we see $H^2(M) = 0$ so $\dim_k H^0(M) \geq 2$ as desired. 

\begin{thm}\label{tmoduliiscurve}
The coarse moduli scheme $\mathcal{M}$ of $A$-line bundles of Chern class $c_1 = L_{\s}\otimes_Y \oy(H), c_2 = 1$ is a smooth genus 2 curve. 
\end{thm}
\textbf{Proof.}
We know from \cite[theorem~2.4]{HoSt} that the coarse moduli scheme $\mathcal{M}$ of torsion-free $A$-modules with the above Chern classes is a projective scheme. Furthermore, proposition~\ref{pwA} ensures that $\mathcal{M}$ is smooth. Note by minimality of the second Chern class, the members of $\mathcal{M}$ are in fact $A$-line bundles. Let $\mathcal{M}_Y$ denote the coarse moduli scheme of (Gieseker) semistable rank two bundles on $Y$ with Chern classes $c_1 = L_{\s}\otimes_Y \oy(H), c_2 = 1$. We certainly have a map $\phi:\mathcal{M} \lm \mathcal{M}_Y$. Also, the $\PP^1$-family of vector bundles $M(p)$ constructed in lemma~\ref{lP1family} gives a map $\psi: \PP^1 \lm \mathcal{M}_Y$ which is injective. Now proposition~\ref{pnis1minc2} shows that $\im \phi = \im \psi$ so functoriality of normalisation ensures a map $\rho:\mathcal{M} \lm \PP^1$. We know $\rho$ is a double cover of $\PP^1$ ramified at 6 points corresponding to the $Y$-split bundles. As $\mathcal{M}$ is smooth, it follows that it is in fact a genus two curve by Riemann-Hurwitz.

\vspace{4mm}
\textbf{Acknowledgements:} Daniel Chan was partially funded by an ARC Discovery project grant.

\end{document}